\documentclass[12pt,epsfig]{article} 
\usepackage{amsbsy,amssymb,amscd,amsfonts,latexsym,amstext,delarray, 
amsmath,epsfig}  \headsep=15pt 
 
\input xypic 
 
\newtheorem{thm}{Theorem}[section] 
\newtheorem{prop}[thm]{Proposition} 
\newtheorem{cor}[thm]{Corollary} 
\newtheorem{lem}[thm]{Lemma} 
 
\newtheorem{defn}[thm]{Definition} 
\newtheorem{rem}[thm]{Remark}

\numberwithin{equation}{section}

\def\R{{\mathbb R}}

\def\m{{\mathfrak m}}

\def\O{{\mathcal O}} 
\def\mX{{\mathfrak X}} 
\def\Sp{{\rm{Spec}}}

\def\Q{{\mathbb Q}} 
\def\C{{\mathbb C}} 
\def\Z{{\mathbb Z}} 
\def\N{{\mathbb N}} 
 
\def\P{{\mathbb P}}

\def\PGL{{\rm PGL}} 
\def\SL{{\rm SL}} 
 
\def\Tr{{\rm Tr}}

\def\mK{{\mathfrak K}} 
 
\def\cW{{\mathcal W}} 
\def\cP{{\mathcal P}} 
\def\cK{{\mathcal K}} 
\def\cS{{\mathcal S}} 
 
\def\cH{{\mathcal H}} 
 
\def\mX{{\mathfrak X}} 
\def\mF{{\mathbb F}} 
\def\cA{{\mathcal A}} 
\def\cV{{\mathcal V}} 
 
\def\mK{{\mathbb K}} 
 
\newcommand{\ie}{{\it i.e.\/}\ } 
\newcommand{\eg}{{\it e.g.\/}\ } 
\newcommand{\cf}{{\it cf.\/}\ }

\title{Spectral triples from Mumford curves} 
\author{C.Consani and M.Marcolli} 
\date{} 
 
\begin{document} 
 
\maketitle 
 
\begin{abstract} 
\noindent We construct spectral triples associated to 
Schottky--Mumford curves, in such a way that the local Euler 
factor can be recovered from the zeta functions of such spectral 
triples. We propose a way of extending this construction to the case
where the curve is not $k$-split degenerate.  
\end{abstract} 
 
\section{Introduction} 
 
Let $X$ be a curve over a finite extension $K$ of $\Q_p$, which is 
$k$-split degenerate, for $k$ the residue field. It is well known 
that such curve admits a $p$-adic uniformization by a $p$-adic 
Schottky group $\Gamma$ acting on the Bruhat-Tits tree $\Delta_K$. 
We associate ${\rm C}^*$-algebras to certain subgraphs $\Delta$ of 
the Bruhat-Tits tree and construct corresponding dynamical 
cohomologies $\cH^1_{dyn}(\Delta/\Gamma)$ that resemble the 
construction at arithmetic infinity given in \cite{CM}. We 
introduce a Dirac operator $D$, which depends on the graded 
structure of the dynamical cohomology, in such a way that the data 
$$ ({\rm C}^*(\Delta/\Gamma), \cH^1_{dyn}(\Delta/\Gamma) \oplus 
\cH^1_{dyn}(\Delta/\Gamma), D) $$ give a spectral triple in the 
sense of Connes. 
 
We recover, from a certain family of zeta functions associated to 
the spectral triple, the local Euler factor $L_p(H^1(X),s)$ of the 
curve $X$, in the form of a regularized determinant as computed by 
Deninger. The advantage of this construction is that it provides, 
in the case we are considering, a natural geometric 
interpretation, in terms of dynamics of walks on the graph 
$\Delta/\Gamma$, of the infinite dimensional cohomology theory 
introduced by Deninger. 

We propose a way of extending the construction to the case where the
curve is not $k$-split degenerate, by enlarging the dual graph of the
special fiber by new edges, in such a way that we also obtain
embeddings in the dynamical cohomology of the first cohomology group
of the components of the special fiber.
 
\section{Directed graphs} 
 
We begin by recalling some generalities about graphs that we will 
use throughout the paper. 
 
\bigskip 
 
 A {\em directed graph} $E$ consists of data $E=(E^0, E^1, E^1_+, 
r, s, \iota)$, where $E^0$ is the set of vertices, $E^1$ is the 
set of oriented edges $w=\{ e, \epsilon \}$, where $e$ is an edge 
of the graph and $\epsilon =\pm 1$ is a choice of orientation. The 
set $E^1_+$ consists of a choice of orientation for each edge, 
namely one element in each pair $\{ e, \pm 1 \}$. The maps $r,s: 
E^1 \to E^0$ are the range and source maps, and $\iota$ is the 
involution on $E^1$ defined by $\iota(w)=\{ e, -\epsilon \}$. 
 
\medskip 
 
A morphism $f$ of directed graphs $E$ and $\tilde E$ consists of 
maps $f^0 : E^0 \to \tilde E^0$ and $f^1: E^1 \to \tilde E^1$ 
which satisfy $f^i \circ r= \tilde r \circ f^i$, $f^i \circ s= 
\tilde s \circ f^i$, and $f^i \circ \iota= \tilde \iota \circ 
f^i$, for $i=0,1$. The morphism $f$ is a monomorphism if the $f^i$ 
are invertible, for $i=0,1$. This defines the automorphism group 
of a directed graph, which we denote $Aut(E)$. A morphism $f$ of 
directed graphs is a covering map if the $f^i$ are onto, for 
$i=0,1$, and $f^1$ gives a bijection $f^1 : s^{-1}(v) 
\stackrel{\sim}{\to} \tilde s^{-1}(f(v))$ and the same with 
respect to the range map. 
 
\medskip 
 
A directed graph is finite if $E^0$ and $E^1$ are finite sets. It 
is {\em row finite} if at each vertex $v\in E^0$ there are at most 
finitely many edges $w$ in $E^1_+$ such that $s(w)=v$. The graph 
is {\em locally finite} if each vertex emits and receives at most 
finitely many oriented edges in $E^1$. A vertex $v$ in a directed 
graph is a {\em sink} if there is no edge in $E^1_+$ with source 
$v$. We denote by $\sigma(E)$ the subset of $E^0$ given by the 
sinks. 
 
\medskip 
 
A juxtaposition of oriented edges $w_1 w_2$ is said to be {\em 
admissible} if $w_2\neq \iota(w_1)$ and $r(w_1)=s(w_2)$. A 
(finite, infinite, doubly infinite) {\em walk} in a directed graph 
$E$ is an admissible (finite, infinite, doubly infinite) sequence 
of elements in $E^1$. A (finite, infinite, doubly infinite) {\em 
path} in $E$ is a walk where all edges in the sequence are in 
$E^1_+$. We denote by $\cW^n(E)$ the set of walks of length $n$, 
by $\cW^*(E)=\cup_n \cW^n(E)$, by $\cW^+(E)$ the set of infinite 
walks, and by $\cW(E)$ the set of doubly infinite walks. 
Similarly, we introduce the analogous notation $\cP^n(E)$, 
$\cP^*(E)$, and $\cP^+(E)$ for paths. We shall drop the explicit 
mention of the graph $E$ in the notation for walks and paths, when 
no confusion arises. We denote by $\sigma^*(E)\subset \cP^*(E)$ 
the set of paths that end at a sink. A directed graph is a 
directed tree if, for any two vertices, there exists a unique walk 
in $\cW^*(E)$ connecting them. 
 
\medskip 
 
The universal cover $\Delta$ of a connected directed graph $E$, 
endowed with a choice of a base point $v_0\in E^0$, is a directed 
tree obtained by setting $\Delta^0=\cW^*(v_0)$, the set of all 
walks in $E$ that start at $v_0$, $\Delta^1=\{ (\omega, w) \in 
\cW^*(v_0) \times E^1 : \, r(\omega)=s(w) \}$, 
$r(\omega,w)=\omega$, $s(\omega,w)=\omega w$, 
$\iota(\omega,w)=(\omega w,\iota(w))$. The fundamental group of 
$E$, with respect to the choice of base point $v_0$, is given by 
$\Gamma =\{ \gamma \in \cW^*(v_0) : \, r(\gamma)=v_0 \}$. 
 
\medskip 
 
Let $G$ be a subgroup of $Aut(E)$. We can then form the quotient 
$E/G$, which is also a directed graph. In particular, if $\Delta$ 
is the universal cover of a directed graph $E$ and $\Gamma$ is the 
fundamental group, with respect to the choice of a base point 
$v_0\in E^0$, then we have an isomorphism of directed graphs 
$E\simeq \Delta/\Gamma$. The map $\Delta \to E$ is a covering map. 
 
\medskip 
 
Two paths $\omega,\tilde \omega$ in $\cP^+$ are {\em shift-tail 
equivalent}, if there exists a $N\geq 1$ and a $k\in \Z$ such that 
$\omega_i=\tilde\omega_{i-k}$ for all $i\geq N$. The boundary of a 
directed tree $\Delta$ is given by $\partial \Delta = (\cP^+ \cup 
\sigma^*)/\sim$, where the shift-tail equivalence is extended to 
the set of paths ending at a sink by the condition $\omega \sim 
\tilde \omega$, for $\omega, \tilde\omega \in \sigma^*$, if and 
only if $r(\omega)=r(\tilde\omega)$. This definition extends to a 
functorial notion of boundary of directed graphs, as shown in 
\cite{Spi}.

\section{Schottky-Mumford curves} 
 
Throughout this chapter we will denote by $p$ a fixed prime number 
and by $\Q_p$ the field of $p$-adic numbers. The field $K$ will be 
a given finite extension of $\Q_p$, with $\O\subset K$ its ring of 
integers, $\m\subset \O$ the maximal ideal and $\pi\in\m$ a 
uniformizer (\ie $\m = (\pi)$). Finally, we will denote by $k$ the 
residue field $k = \O/\m$. This is a finite field of cardinality 
$q = \rm{card}(\O/\m)$. 
 
Let $V$ be a two-dimensional vector space over $K$. We write $\P^1(K)$ 
for the set of $K$-rational points of $\P^1_K$, the projective line 
over $K$. The space $\P^1(K)$ is identified with the set of $K$-lines 
passing through the origin in $V$. Let $G = \PGL(2,K)$ be the group of 
automorphisms of $\P^1(K)$. 
 
In this first chapter we collect some well known facts and properties 
on the tree of the group $\PGL(2,K)$ and on the action of a Schottky 
group on such tree. Detailed explanations for the statements we will 
recall here without proofs are contained in \cite{Ma} and \cite{Mu}. 
 
\subsection{The tree of the group \PGL(2,K)} \label{treeSec} 
 
The description of the vertices of the tree associated to $\PGL(2,K)$ 
is as follows. One considers the set of free $\O$-modules of rank $2$: 
$M\subset V$. Two such modules are {\it equivalent} $M_1 \sim M_2$ if 
there exists an element $\lambda\in K^\ast$, such that $M_1 = \lambda 
M_2$. 
The group $\rm{GL}(V)$ of linear automorphisms of $V$ operates on the 
set of such modules {\it on the left}: $gM = \{gm~|~m\in M\}$, $g\in 
\rm{GL}(V)$. Notice that the relation $M_1 \sim M_2$ is equivalent to 
the condition that $M_1$ and $M_2$ belong to the same orbit of the center 
$K^\ast \subset \rm{GL}(V)$. Hence, the group $G = \rm{GL}(V)/K^\ast$ 
operates (on the left) on the set of classes of equivalent modules. 
 
We denote by $\Delta^0_K$ the set of such classes and by $\{M\}$ 
the class of the module $M$. Because $\O$ is a principal ideals 
domain and every module $M$ has two generators, it follows that 
\[ 
\{M_1\}, \{M_2\} \in\Delta^0_K, M_1\supset 
M_2\quad\Rightarrow\quad M_1/M_2 \simeq \O/\m^l \oplus 
\O/\m^k,\quad l,k\in\N. 
\] 
 
The multiplication of $M_1$ and $M_2$ by elements of $K$ preserves the 
inclusion $M_1\supset M_2$, hence the natural number 
\begin{equation}\label{dist} 
d(\{M_1\},\{M_2\}) = \vert l-k\vert 
\end{equation} 
is well defined. 
 
\begin{defn} The graph $\Delta_K$ of the group $\PGL(2,K)$ is the 
infinite graph with set of vertices $\Delta^0_K$, in which two 
vertices $\{M_1\},\{M_2\}$ are adjacent and hence connected by an 
edge if and only if $d(\{M_1\},\{M_2\}) = 1$. 
\end{defn} 
 
The following properties characterize completely $\Delta_K$ and 
are well known (\cf~\eg\cite{Ma} and \cite{Mu}). 
 
\begin{prop} \begin{enumerate} 
\item The graph $\Delta_K$ is a connected, locally finite tree with  $q+1$ 
edges departing from each of its vertices. 
\item The shortest walk in $\Delta_K$ connecting two vertices $\{M_1\}, 
\{M_2\}$ of $\Delta^0_K$ without back-tracking has length 
$d(\{M_1\},\{M_2\})$. 
\item The group $G$ acts (on the left) transitively on $\Delta^0_K$ and 
it preserves the metric $d$. 
\end{enumerate} 
\end{prop} 
 
The tree $\Delta_K$ is called the Bruhat--Tits tree associated to 
the group $G = \PGL(2,K)$.\vspace{.1in} 
 
A {\it half-line} in $\Delta_K$ is an infinite sequence 
$(\{M_n\})_{n\in\N}$ of vertices of $\Delta_K$ {\it without 
repetitions} such that $\{M_n\}$ is adjacent to $\{M_{n-1}\}$ for 
all $n$. Thus, a half-line is given by a sequence $M_0\supset 
M_1\supset\ldots$ of free $\O$-modules where $M_0/M_n \simeq 
\O/(\pi^n)$ for all $n$. 
 
The subspace $K(\cap_{n\in\N}M_n)\subset V$ defines a point of 
$\P^1(K)$. Conversely, given a point of $\P^1(K)$ represented by a 
vector $v_1\in V$, choose a second vector $v_2\in V$ such that 
$\{v_1,v_2\}$ form a basis of $V$. Let $M_n$ be the free 
$\O$-module with basis $\{v_1,\pi^nv_2\}$. Then, 
$K(\cap_{n\in\N}M_n)$ defines the point of $\P^1(K)$ we started 
with. 
 
Two half-lines are said to be equivalent if they differ only by a 
finite number of vertices. An equivalence class of halflines is 
called an {\it end} of $\Delta_K$. The set of ends of $\Delta_K$ 
will be denoted by $\partial\Delta_K$ (the ``boundary'' of 
$\Delta_K$). We shall give $\Delta_K$ the structure of a directed 
graph, in such a way that this notion of boundary agrees with the 
one described in the previous section. 
 
It is immediate to verify that the construction described above 
defines a one-to-one correspondence 
$\partial\Delta_K\quad\rightleftarrows\quad\P^1(K)$ between 
equivalence classes of half-lines 
and elements of $\P^1(K)$. 
 
\subsection{The action of a Schottky group on the tree $\Delta_K$} 
 
A Schottky group $\Gamma$ is a subgroup of $\PGL(2,K)$ which is 
finitely generated and whose elements $\gamma\neq 1$ are {\it 
hyperbolic} (\ie the eigenvalues of $\gamma$ in $K$ have different 
valuation). The group $\Gamma$ is discrete in $G = \PGL(2,K)$, 
torsion free and acts freely on the tree $\Delta_K$. Furthermore, 
one can show that $\Gamma$ acts discretely at some point $z\in 
\P^1(\bar K)$. 
 
Let $\Lambda_\Gamma\subset \P^1(K)$ be the closure of the set of 
points in $\P^1(K)$ that are fixed points of some 
$\gamma\in\Gamma\setminus\{1\}$. This is called the {\em limit 
set} of $\Gamma$. We have $\rm{card}(\Lambda_\Gamma) < \infty$ if 
and only if $\Gamma = (\gamma)^\Z$, for some $\gamma\in\Gamma$. We 
will denote by $\Omega_\Gamma=\Omega_\Gamma(K)$ the set of points 
on which $\Gamma$ acts discretely; equivalently said 
$\Omega_\Gamma$ is the set of points which are not limits of fixed 
points of elements of $\Gamma$: $\Omega_\Gamma = 
\P^1(K)\setminus\Lambda_\Gamma$. This set is called the {\it 
domain of discontinuity} for the Schottky group $\Gamma$. 
 
A path in $\Delta_K$, infinite in both directions and with no 
back-tracking, is called an {\it axis} of $\Delta_K$. Any two 
points $z_1,z_2\in \P^1(K)$ uniquely define their connecting axis 
whose endpoints lie in the classes described by $z_1$ and $z_2$ in 
$\partial\Delta_K$. 
 
Any hyperbolic element $\gamma\in\Gamma$ has two fixed points in 
$\P^1(K)$. The unique axis of $\Delta_K$ whose ends are fixed by 
$\gamma$ is called the axis of $\gamma$. The element $\gamma$ acts 
on its axis as a translation. 
 
\begin{figure} 
\begin{center} 
\epsfig{file=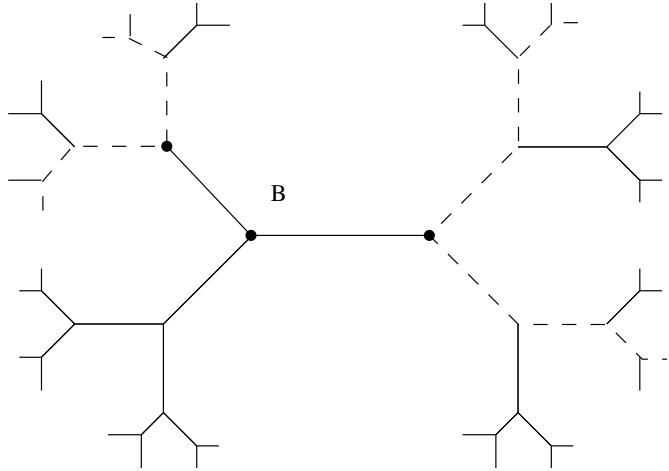} 
\caption{Two axes and the bridge $B$ between them. \label{axes}} 
\end{center} 
\end{figure}

Suppose two axes are given in $\Delta_K$. Any path without 
back-tracking beginning on one axis and ending on the other and 
with no edges in common with either axis is said to be the {\it 
bridge} between them. A bridge may consist of a single point, else 
it is uniquely defined (Figure \ref{axes}). 
 
For any Schottky group $\Gamma \subset \PGL(2,K)$ there is a 
smallest {\it subtree} $\Delta'_\Gamma \subset \Delta$ containing 
the axes of all elements of $\Gamma$. Equivalently said, 
$\Delta'_\Gamma$ is the maximal connected subgraph of $\Delta_K$ 
containing the axes of all elements of $\Gamma$ and the bridges 
between them. 
 
The set of ends of $\Delta'_\Gamma$ in $\P^1(K)$ is 
$\Lambda_\Gamma$, the limit set of $\Gamma$. The group $\Gamma$ 
carries $\Delta'_\Gamma$ into itself so that the quotient 
$\Delta'_\Gamma /\Gamma$ is a {\it finite graph}. 
 
The graph $\Delta'_\Gamma /\Gamma$ has an important geometric 
interpretation as the dual graph of the closed fiber of the {\it 
minimal smooth model} over $\O$ ($k$-split degenerate semi-stable 
curve) of the algebraic curve $C$ holomorphically isomorphic to 
$X_\Gamma := \Omega_\Gamma /\Gamma$ (\cf~\cite{Mu} p.~163). 
 
Furthermore, there is a smallest tree $\Delta_\Gamma$ on which 
$\Gamma$ acts, with vertices $\Delta_\Gamma^{0}\subset 
\Delta^0_K$, and such that the finite graph $\Delta_\Gamma/ 
\Gamma$ is the dual graph of the specialization of $C$ over $\O$. 
The set $\Delta_\Gamma^{0}$ is a subset of the set of vertices of 
$\Delta'_\Gamma$. 
 
The degenerating curve $C$ describing the analytic uniformization 
$X_\Gamma \stackrel{\simeq}\to C$ is a {\it $k$-split degenerate, 
stable curve}. When the genus of the fibers is at least 2 - \ie when 
the Schottky group has at least $g\geq 2$ generators - the curves 
$X_\Gamma$ are called Schottky--Mumford curves. 
 
\subsection{Field extensions: functoriality} 
 
Let $L\supset K$ be a finite fields extension with ramification index 
$e_{L/K}$ and rings of integers $\O_L$ and $\O_K$. If $M \subset V$ is 
a free $\O_K$-module of rank $2$, then $M\otimes_{\O_K}\O_L \subset 
V\otimes_K L$ is a free $\O_L$-module of the same rank. It is obvious 
that equivalent modules remain equivalent, hence one gets a natural 
embedding of the sets of vertices $\Delta_K^0 \hookrightarrow 
\Delta_L^0$. 
 
The isomorphism $(\O_K/\m^r)\otimes \O_L \simeq \O_L/\m^{re_{L/K}}$ 
shows that distance would not, in general, be preserved 
under field extensions. To eliminate this disadvantage, one introduces 
on graphs $\Delta_L$, for all extensions $L\supset K$, a 
$K$-normalized distance 
\[ 
d_K(\{M_1\},\{M_2\}) = \frac{1}{e_{L/K}}d_L(\{M_1\},\{M_2\});\qquad 
M_1, M_2 \subset V\otimes L. 
\] 
 
\begin{figure} 
\begin{center} 
\epsfig{file=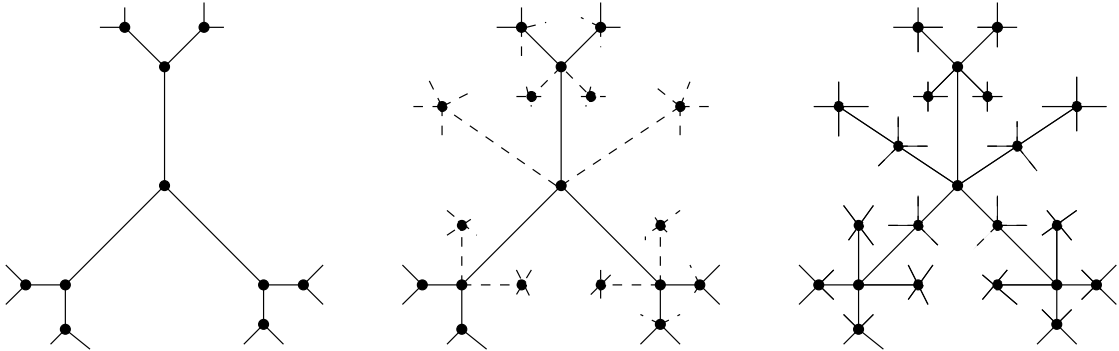} 
\caption{The tree $\Delta_K$ for $K=\Q_2$ and $\Delta_L$ for field 
extensions with $f=2$ and $e_{L/K}=2$ \label{extTree}} 
\end{center} 
\end{figure}

This way, the embedding $\Delta_K^0 \hookrightarrow \Delta_L^0$ 
becomes isometric. When $L \supset K$ ramifies, $e_{L/K}-1$ new 
vertices appear in $\Delta^0_L$ between each couple of adjacent 
vertices in $\Delta^0_K$ (\cf the third graph in Figure 
\ref{extTree}). Moreover, $q^f + 1$ edges each of which of 
length $\frac{1}{e_{L/K}}$, for $f = \frac{1}{e_{L/K}}[L:K]$, depart 
from each vertex in $\Delta^0_L$ (\cf the second graph in Figure 
\ref{extTree}). 
 
Because $\PGL(2,K)\subset \PGL(2,L)$ and the natural embedding 
$\P^1(K) \subset \P^1(L)$ is compatible with a concept of 
$K$-direction of exiting from any vertex of $\Delta^0_K$, the 
construction is functorial under finite fields extensions and this 
process determines, for a fixed Schottky group $\Gamma \subset 
\PGL(2,K)$, a projective system $\{X_{L,\Gamma}~:~[L:K]<\infty\}$ of 
Schottky-Mumford curves. 
 
\subsection{Edge orientation} 
 
\begin{figure} 
\begin{center} 
\epsfig{file=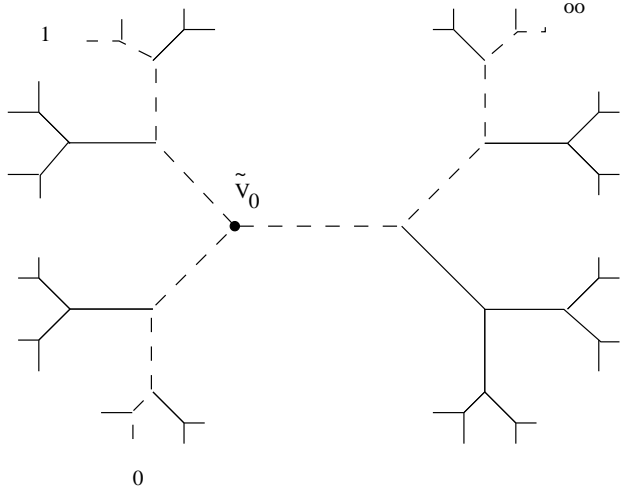} \caption{The crossroad $\tilde v_0$ of the 
points $\{ 0, 1, \infty \}$ \label{cross}} 
\end{center} 
\end{figure}

We now show how to endow the graphs $\Delta_\Gamma$ and $\Delta_K$ 
with the structure of directed graph. The choice of a coordinate 
$z\in \P^1(K)$ determines a base point $\tilde v_0$ in $\Delta_K$. 
In fact, the points $\{ 0, 1, \infty \}$ on $\P^1(K)$ determine a 
unique {\em crossroad}: the unique vertex of $\Delta_K$ with the 
property that the walks without back-tracking from $\tilde v_0$ to 
the three points on $\P^1(K)$ start from $\tilde v_0$ in three 
different directions (Figure \ref{cross}).

In order to obtain a structure of directed graph on the tree 
$\Delta_K$, for each $z\in \P^1(K)$ we consider the unique 
infinite chain of edges without backtracking in $\Delta_K$ that 
has initial vertex $\tilde v_0$ and whose equivalence class modulo 
shift--tail equivalence is the point $z$. We declare such chain of 
edges to be a path in $\cP^+(\Delta_K)$. This determines on 
$\Delta_K$ the structure of a directed graph, with $\partial 
\Delta_K = \P^1(K)$. This agrees with the boundary as defined in 
\S \ref{treeSec}. 
 
We assume here that the coordinate $z\in \P^1(K)$ is chosen in 
such a way that the crossroad $\tilde v_0$ of $\{ 0, 1, \infty \}$ 
is a vertex of $\Delta_\Gamma$. Then, by the same procedure, 
we can regard $\Delta_\Gamma$ as a directed graph with $\partial 
\Delta_\Gamma =\Lambda_\Gamma$.

\section{$C^*$-algebras of graphs}\label{Agraph} 
 
In this section we recall the construction of ${\rm C}^*$-algebras 
associated to locally finite directed graphs. We follow mostly the 
references \cite{BPRS}, \cite{KuPa}, \cite{KuPaRae}, \cite{Spi}. 
We refer the reader to the bibliography of the 
aforementioned articles for more information. For simplicity, we 
state the following results in the special case of locally finite 
directed graph, though the theory extends to more general directed 
graphs (\cf \eg \cite{Spi}). 
 
\bigskip 
 
 A {\em Cuntz--Krieger family} consists of a collection $\{ P_v 
\}_{v\in E^0}$ of mutually orthogonal projections and $\{
S_w \}_{w \in E^1_+ }$ of partial isometries, satisfying the 
conditions: $S_w^* S_w =P_{r(w)}$ and, for all $v\in s(E^1_+)$, 
$P_v = \sum_{w: s(w)=v} S_w S_w^*$. 
 
\medskip 
 
The {\em edge matrix} $A_+$ of a locally finite (or row finite) 
directed graph is an $(\# E^1_+)\times (\# E^1_+)$ (possibly 
infinite) matrix. The entries $A_+(w_i,w_j)$ satisfy $A_+(w_i,w_j)=1$ 
if $w_i w_j$ is an admissible path, and $A_+(w_i,w_j)=0$ otherwise. 
The Cuntz--Krieger elements $\{ P_v, S_w \}$ satisfy the relation 
$S_w^* S_w = \sum A(w,\tilde w) S_{\tilde w} S_{\tilde w}^*$. The 
{\em directed edge matrix} of $E$ is a $\# E^1 \times \# E^1$ 
(possibly infinite) matrix with entries $A(w_i,w_j)=1$ if $w_i 
w_j$ is an admissible walk and $A(w_i,w_j)=0$ otherwise. 
 
\medskip 
 
There is a universal $C^*$-algebra $C^*(E)$ generated by a 
Cuntz--Krieger family. If $E$ is a finite graph with no sinks, we 
have $C^*(E)\simeq {\mathcal O}_{A_+}$, where ${\mathcal O}_{A_+}$ 
is the Cuntz-Krieger algebra of the edge matrix $A_+$. If the 
directed graph is a tree $\Delta$, then $C^*(\Delta)$ is an  AF 
algebra strongly Morita equivalent to the commutative 
$C^*$-algebra $C_0(\partial \Delta)$. A monomorphism of directed 
trees induces an injective $*$-homomorphism of the corresponding 
$C^*$-algebras. 
 
\medskip 
 
If $G\subset Aut(E)$ is a group acting freely on the directed 
graph $E$, with quotient graph $E/G$, then the crossed product 
$C^*$-algebra $C^*(E)\rtimes G$ is strongly Morita equivalent to 
$C^*(E/G)$. In particular, if $\Delta$ is the universal covering 
tree of a directed graph $E$ and $\Gamma$ is the fundamental 
group, then the algebra $C^*(E)$ is strongly Morita equivalent to 
$C_0(\partial \Delta)\rtimes \Gamma$. 
 
\medskip 
 
There is a {\em gauge action} of $U(1)$ on the graph algebra 
$C^*(E)$ given by $\lambda : \{ P_v, S_w \} \mapsto \{ P_v, 
\lambda S_w \}$. A subset $H$ of the set of vertices $E^0$ of a 
directed graph is {\em saturated hereditary} if $v\in H$ implies 
that, for all $\omega \in \cP^*(E)$ with $s(\omega)=v$, also 
$r(\omega)\in H$, and conversely, if any $\omega \in \cP^*(E)$ 
with $s(\omega)=v$ satisfies $r(\omega)\in H$, then also $v\in H$. 
For a locally finite graph there is a bijective correspondence 
between saturated hereditary subsets of $E^0$ and gauge invariant 
closed ideals in $C^*(E)$. In the case of a tree $\Delta$, there 
is a bijection between saturated hereditary subsets of $\Delta^0$ 
and open sets in $\partial \Delta$. This is proved for the more 
general (non necessarily locally finite) case in \cite{Spi}. 
 
\medskip 
 
It is convenient to consider also the Toeplitz extensions (\cf 
\cite{Spi}) 
$$ 0\to I_S \to T\O(E,S) \to C^*(E) \to 0, $$ 
where $S$ is a subset of $E^0$ and $I_S= \oplus_{v\in S^c} \cK_v$, 
where $\cK_v$ is the algebra of compact operators on a Hilbert 
space of dimension $\#(\cP^* \cap r^{-1}(v))$. The $C^*$-algebra 
$T\O(E,S)$ is generated by operators $\{ S_w \}_{w\in E^1_+}$ and 
$\{ P_v \}_{v\in E^0}$ satisfying the conditions: $S_w^* S_w 
=P_{r(w)}$ and, for all $v\in s(E^1_+)$, $P_v \geq \sum_{w: 
s(w)=v} S_w S_w^*$, with equality for $v\in S$. 
 
If $j: E \hookrightarrow \tilde E$ is an inclusion of directed 
graphs, the following functoriality property holds (\cf 
\cite{Spi}): given $\tilde S \subset \tilde E$, consider the set 
$S$ of vertices $v$ in $E^0$ such that $j(v)\in \tilde S$ and the outgoing 
edges in $\tilde E^1_+$ with origin at $j(v)$ are all of the form 
$j(w)$, for some $w\in E^1_+$ with $s(w)=v$. Then this induces an 
injective $*$-homomorphism $J : T\O(E,S) \to T\O(\tilde E, \tilde 
S)$. 
 
In particular, for a family of subgraphs $E$ of a directed graph 
$\tilde E$, ordered by inclusions, with $\cup E^0= \tilde E^0$ and 
$\cup E^1 = \tilde E^1$, and for a choice of $\tilde S\subset 
\tilde E^0$ and corresponding $S$ as above, we have 
\begin{equation}\label{dir-lim} 
T\O(\tilde E, \tilde S) = \lim_{E} T\O(E,S). 
\end{equation}

\subsection{Reduction graphs} 
 
In the theory of Mumford curves, it is important to consider also 
the reduction modulo powers $\m^n$ of the maximal ideal $\m 
\subset \O_K$, which provides infinitesimal neighborhoods of 
order $n$ of the closed fiber. In the language of $C^*$-algebras, 
this corresponds to the following construction. 
 
For each $n\geq 0$, we consider a subgraph $\Delta_{K,n}$ of the 
Bruhat-Tits tree $\Delta_K$ defined by setting 
$$ \Delta_{K,n}^0 := \{ v\in \Delta_K^0 : \,  d(v, 
\Delta_\Gamma ')\leq n \}, $$ with respect to the distance 
\eqref{dist}, with $d(v,\Delta_\Gamma ') := \inf \{ d(v,\tilde v): 
\, \tilde v \in (\Delta_\Gamma ')^0 \}$, and 
$$ \Delta_{K,n}^1 :=\{ w\in \Delta_K^1 : \, s(w), r(w)\in 
\Delta_{K,n}^0 \}. $$ Thus, we have $\Delta_{K,0} = \Delta_\Gamma 
'$. For $n\geq 1$ the graph $\Delta_{K,n}$ has a non-empty set of 
sinks $\sigma_{K,n}\subset\partial \Delta_{K,n}$. We have 
$\Delta_K=\cup_n \Delta_{K,n}$. 
 
For all $n\in \N$, the graph $\Delta_{K,n}$ is invariant under the 
action of the Schottky group $\Gamma$ on $\Delta$, and the finite 
graph $\Delta_{K,n}/\Gamma$ gives the dual graph of the reduction 
$X_K \otimes \O/ \m^{n+1}$. Thus, we refer to the $\Delta_{K,n}$ 
as {\em reduction graphs}. They form a directed family with 
inclusions $j_{n,m}: \Delta_{K,n} \hookrightarrow \Delta_{K,m}$, 
for all $m\geq n$, with all the inclusions compatible with the 
action of $\Gamma$. 
 
For each reduction graph, we can consider corresponding ${\rm 
C}^*$-algebras ${\rm C}^*(\Delta_{K,n})$ and ${\rm 
C}^*(\Delta_{K,n}/\Gamma)\simeq {\rm C}^*(\Delta_{K,n})\rtimes 
\Gamma$ (Morita equivalence). 
 
The following result, which is a direct application of the 
statements on graph ${\rm C}^*$-algebras listed in \S 
\ref{Agraph}, describes the relation between the algebras ${\rm 
C}^*(\Delta_{K,n}/\Gamma)$ and ${\rm C}^*(\Delta_K/\Gamma)$. 
 
\begin{lem} 
We have injective $*$-homomorphisms $J_{n,m}: C^*(\Delta_{K,n}) 
\to C^*(\Delta_{K,m})$. Correspondingly, if we set 
$S_n=\Delta_{K,n}^0 \backslash \sigma_{K,n}$, we obtain 
$$ C^*(\Delta_K/\Gamma) = \lim_n T\O(\Delta_{K,n}/\Gamma, 
S_n/\Gamma), $$ 
where $T\O(\Delta_{K,n}/\Gamma,S_n/\Gamma)$ 
satisfies 
$$ 0 \to \oplus_{v\in \sigma(\Delta_{K,n}/\Gamma)} \cK_v \to 
T\O(\Delta_{K,n}/\Gamma,S_n/\Gamma) \to C^*(\Delta_{K,n}/\Gamma) 
\to 0. $$ 
\end{lem}

\section{Dynamics of walks on dual graphs} 
 
In this section we introduce a dynamical system associated to the 
space $\cW(\Delta/\Gamma)$ of walks on the quotient of a directed tree
$\Delta$ by a free action of $\Gamma$. In particular, we are
interested in the  
cases where $\Delta$ is (a certain extension of) one of the graphs 
$\Delta_{K,n}$ for some $n\geq 0$. The same construction applies 
to the tree $\Delta_\Gamma$, where this dynamical system is a 
subshift of finite type associated to the action of the Schottky 
group $\Gamma$ on its limit set $\Lambda_\Gamma$, analogous to the 
one considered in \cite{CM}. 
 
\medskip 
 
Let $\bar V \subset \Delta$ be a finite subtree whose set of edges 
consists of one representative for each $\Gamma$-class. This is a 
fundamental domain for $\Gamma$ in the weak sense (following the 
notation of \cite{Ma}), since some vertices may be identified 
under the action of $\Gamma$. Correspondingly, we denote by $V$ 
the set of ends of all infinite paths in $\Delta$ starting at 
points in $\bar V$. 
 
We assume that, for the $\Gamma$-invariant directed tree $\Delta$, 
the set $\bar V$ has {\em finitely many edges}. This is the case 
for $\Delta_\Gamma$ as well as for any of the trees 
$\Delta_{K,n}$. 
 
Since, for $n\geq 1$, the graphs $\Delta_{K,n}$ have sinks, in 
order to consider the space of doubly infinite walks on these 
graphs, we need to complete each walk ending at a sink to an 
infinite walk obtained by repeating the last word. This is 
equivalent to extending the graph $\Delta_{K,n}$ by adding an 
infinite tail to each sink. Appending tails to sinks is standard 
technique in the theory of graph ${\rm C}^*$-algebras, in order to 
reduce the general case to the easier case of graphs with no 
sinks. We use the notation $\bar \Delta_{K,n}$ for the completed 
graph with infinite tails. These have an action of $\Gamma$ 
obtained by extending the action on $\Delta_{K,n}$, by translating 
the whole tail in the same way as the corresponding sink in 
$\Delta_{K,n}$. 
 
For $\Delta=\bar\Delta_{K,n}$, the set $\cW(\Delta/\Gamma)$ of 
doubly infinite walks on the graph 
$\bar\Delta_{K,n}/\Gamma$ can be identified with the set of doubly 
infinite admissible sequences in the {\em finite alphabet} given 
by the edges in the fundamental domain $\bar V$ of the graph 
$\Delta_{K,n}$, with both possible orientations. 
 
On $\cW(\Delta/\Gamma)$ we consider the topology generated by the 
sets $\cW^s(\omega,\ell)=\{ \tilde\omega \in \cW(\Delta/\Gamma): 
\tilde\omega_k=\omega_k, k\geq \ell\}$ and $\cW^u(\omega,\ell)=\{ 
\tilde\omega \in \cW(\Delta/\Gamma): \tilde\omega_k=\omega_k, 
k\leq \ell\}$, for $\omega\in \cW(\Delta/\Gamma)$ and $\ell \in 
\Z$. With this topology, the space $\cW(\Delta/\Gamma)$ is a 
totally disconnected compact Hausdorff space. The invertible shift 
map $T$, given by $(T\omega)_k = \omega_{k+1}$, is a homeomorphism 
of $\cW(\Delta/\Gamma)$. 
 
We have just described the dynamical system 
$(\cW(\Delta/\Gamma),T)$ in terms of subshifts of finite type, 
according to the following definition: 
 
\begin{defn} 
A {\em subshift of finite type} $(\cS_A,T)$ consists of all doubly 
infinite sequences in the elements of a given finite set $W$ 
(alphabet) with the admissibility condition specified by a $\# W 
\times \# W$ elementary matrix, 
$$ \cS_A =\{ (w_k)_{k\in \Z}: w_k \in W, A(w_k,w_{k+1})=1 \}, $$ 
with the action of the invertible shift $(Tw)_k =w_{k+1}$. 
\end{defn} 
 
\begin{lem}\label{subshift} 
The space $\cW(\Delta/\Gamma)$ with the action of the invertible 
shift $T$ is a subshift of finite type, where 
$\cW(\Delta/\Gamma)=\cS_A$ with $A$ the directed edge matrix of 
the finite graph $\Delta/\Gamma$. 
\end{lem} 
 
We can consider the mapping torus of $T$, 
\begin{equation}\label{map-torus-G} 
\cW(\Delta/\Gamma)_{T}:= \cW(\Delta/\Gamma) \times [0,1] / 
(Tx,0)\sim (x, 1). 
\end{equation} 
 
We consider the first cohomology group 
$H^1(\cW(\Delta/\Gamma)_{T},\Z)$, identified with the group of 
homotopy classes of continuous maps of $\cW(\Delta/\Gamma)_{T}$ to 
the circle. If we denote by ${\rm C}(\cW(\Delta/\Gamma),\Z)$ the 
$\Z$-module of integer valued continuous functions on 
$\cW(\Delta/\Gamma)$, and by ${\rm C}(\cW(\Delta/\Gamma),\Z)_T$ 
the cokernel of the map $\delta(f)=f-f\circ T$, we obtain the 
following result (\cf \cite{BoHa}, \cite{PaTu}). 
 
\begin{prop}\label{filtration} 
The map $f\mapsto [\exp(2\pi i t f(x))]$, which associates to an 
element $f\in {\rm C}(\cW(\Delta/\Gamma),\Z)$ a homotopy class of 
maps from $\cW(\Delta/\Gamma)_{T}$ to the circle, gives an 
isomorphism ${\rm C}(\cW(\Delta/\Gamma),\Z)_T \simeq 
H^1(\cW(\Delta/\Gamma)_{T},\Z)$. Moreover, there is a filtration 
of ${\rm C}(\cW(\Delta/\Gamma),\Z)_T$ by free $\Z$-modules 
$F_0\subset F_1\subset \cdots F_n \cdots$, of rank $\theta_n - 
\theta_{n-1} +1 $, where $\theta_n$ is the number of admissible 
words of length $n+1$ in the alphabet, so that we have 
$$ H^1(\cW(\Delta/\Gamma)_{T},\Z)=\varinjlim_n 
F_n. $$ 
The quotients $F_{n+1}/F_n$ are also torsion free. 
\end{prop} 
 
\noindent\underline{Proof.} A continuous function $f\in {\rm 
C}(\cW(\Delta/\Gamma),\Z)$ depends on just finitely many 
coordinates $\omega_k$ of $\omega\in \cW(\Delta/\Gamma)$. In 
particular, this implies that, for some $k_0$, the composite 
$f\circ T^{k_0}$ is a function of only `future coordinates' 
($\omega_k$ with $k\geq 0$). We denote by $\cP\subset {\rm 
C}(\cW(\Delta/\Gamma),\Z)$ the submodule of functions of future 
coordinates. It is then clear that we have ${\rm 
C}(\cW(\Delta/\Gamma),\Z)_T \simeq \cP/\delta \cP$. We also have 
an identification $\cP \simeq {\rm C}(\cW^+(\Delta/\Gamma),\Z)$. 
This gives a filtration $\cP=\cup_n \cP_n$, where $\cP_n$ is 
identified with the submodule of ${\rm 
C}(\cW^+(\Delta/\Gamma),\Z)$ generated by characteristic functions 
of $\cW^+(\Delta/\Gamma,\rho)\subset \cW^+(\Delta/\Gamma)$, where 
$\rho\in \cW^*(\Delta/\Gamma)$ is a finite walk $\rho=w_0\cdots 
w_n$ of length $n+1$, and $\cW^+(\Delta/\Gamma,\rho)$ is the set 
of infinite paths $\omega\in \cW^+(\Delta/\Gamma)$, with $\omega_k 
=w_k$ for $0\leq k\leq n+1$. We have $\delta: \cP_n \to 
\cP_{n+1}$, with kernel the constant functions. We set $F_n := 
\cP_n /\delta(\cP_{n-1})$, for $n\geq 1$ and $F_0 = \cP_0$. The 
inclusions $\cP_n \subset \cP_{n+1}$ induce injections $j : F_n 
\hookrightarrow F_{n+1}$, $j(f+\delta \cP_{n-1})= f+\delta \cP_n$, 
such that $\cP/\delta \cP =\varinjlim_n F_n$. As $\Z$-modules, 
both the $F_n$ and the quotients $F_{n+1}/ j(F_n)$ are torsion 
free. 
 
\noindent $\diamond$ 
 
For more details see the analogous argument given in \cite{PaTu} 
Thm. 19 p. 62-63. 
 
\subsection{The effect of field extensions}

Let $L \supset K$ be a finite extension, with branching index 
$e_{L /K}$ and with $f= [L:K]/e_{L /K}$. Then there is an 
embedding of the set of vertices $\Delta^{(0)}_{K} \subset 
\Delta^{(0)}_{L}$. In between every two vertices of 
$\Delta^{(0)}_{K} \subset \Delta^{(0)}_{L}$ there are $e_{L /K}-1$ 
new vertices of $\Delta_{L}$. If in $\Delta_{K}$ every vertex has 
valence $q+1$, then every vertex of $\Delta_{L}$ has valence $q^f 
+1$. 
 
In particular, the image of the tree $\Delta_\Gamma \subset 
\Delta_{K}$ in $\Delta_{L}$ is the tree $\Delta_\Gamma \subset 
\Delta_{L}$, whereas, when considering the tree $\Delta_\Gamma'$, 
we are inserting $e_{L /K}-1$ new vertices in between each two 
vertices of $\Delta_\Gamma' \subset \Delta_{K}$.  Notice that the 
algebras $C^*(\Delta_\Gamma')$ and $C^*(\Delta_\Gamma)$ are 
strongly Morita equivalent, since they both are equivalent to the 
commutative $C^*$-algebra $C(\Lambda_\Gamma)$. Thus we have: 
 
\begin{lem}{\em 
The strong Morita equivalence class of the graph $C^*$-algebras 
$C^*(\Delta_\Gamma')$ and $C^*(\Delta_\Gamma'/\Gamma)$ is 
independent of finite field extensions $L \supset K$. } 
\end{lem} 
 
\begin{figure} 
\begin{center} 
\epsfig{file=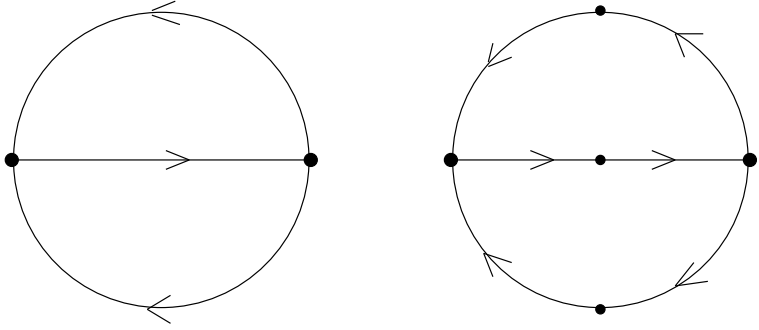} \caption{Effect of a field extension with 
$e_{L /K}=2$. \label{graph2}} 
\end{center} 
\end{figure}

The following example illustrates how the algebra $\O_{A_+}$ 
changes under field extensions, by showing the change in the edge 
matrix $A_+$. If the first directed graph shown in Figure 
\ref{graph2} arises as dual graphs of the closed fiber for a 
totally split degenerate curve, then the effect on the graph of a 
field extension with $e_{L/K}=2$ is illustrated in the second 
graph in Figure \ref{graph2}. 
 
The edge matrix of the original graph was of the form 
$$ A_+= \left( \begin{array}{ccc} 0&1&0\\1&0&1\\0&1&0\end{array} 
\right), $$ while the new edge matrix becomes: 
$$ A_+= \left(\begin{array}{cccccc} 0&1&0&0&0&0\\0&0&1&0&0&0\\ 
0&0&0&1&0&0\\1&0&0&0&1&0\\0&0&0&0&0&1\\0&0&1&0&0&0 \end{array}\right). 
$$ 
 
\bigskip

We are interested in understanding the effect of field extensions 
on the construction considered in the previous paragraph. We have 
the following result. 
 
\begin{prop} 
A finite field extension $L\supset K$ determines a 
morphism 
$$ H^1(\cW(\bar\Delta_{L,n}/\Gamma)_T) \to 
H^1(\cW(\bar\Delta_{K,n}/\Gamma)_T), $$ 
which is compatible with the filtrations. 
\end{prop} 
 
\noindent\underline{Proof.} In \cite{Ma} it is shown that there is 
a canonical choice of the fundamental domains $\bar V$ and $V$ for 
the action of $\Gamma$ on the Bruhat-Tits tree $\Delta_K$ that is 
functorial under a finite extensions $L\supset K$. This determines 
corresponding functorial choices of fundamental domains for the 
graphs $\Delta_{K,n}$. We obtain this way a corresponding 
embedding of walk spaces 
$$ J_{L,K,n}: \cW(\bar\Delta_{K,n}/\Gamma)\hookrightarrow 
\cW(\bar\Delta_{L,n}/\Gamma), $$ 
which replaces each edge in a sequence $\omega \in 
\cW(\bar\Delta_{K,n}/\Gamma)$ (an edge in $\bar V$ for 
$\Delta_{K,n}$) with the corresponding $e_{L/K}$ consecutive edges 
in the fundamental domain $\bar V$ for $\Delta_{L,n}$, thus 
obtaining an element in $\cW(\bar\Delta_{L,n}/\Gamma)$. This map 
satisfies $J_{L,K,n} \circ T = T^{e_{L/K}} \circ J_{L,K,n}$. Thus, 
we obtain an induced map 
$$ 
J_{L,K,n,T}: \cW(\bar\Delta_{K,n}/\Gamma)_T \to 
\cW(\bar\Delta_{L,n}/\Gamma)_{T^{e_{L/K}}}, 
$$ 
where, for $\ell\geq 1$, $\cW(\bar\Delta_{L,n}/\Gamma)_{T^\ell}$
denotes the mapping torus  
$$ \cW(\bar\Delta_{L,n}/\Gamma)_{T^\ell} \simeq \cW(\bar\Delta_{L,n}/\Gamma) 
\times [0,\ell] / (x,0)\sim (T^\ell x, \ell), $$ with a covering 
map $\pi_\ell:\cW(\bar\Delta_{L,n}/\Gamma)_{T^\ell} \to 
\cW(\bar\Delta_{L,n}/\Gamma)_T$. Thus, we obtain a map 
$$ 
\pi_{e_{L/K}} \circ J_{L,K,n,T} : \cW(\bar\Delta_{K,n}/\Gamma)_T 
\to \cW(\bar\Delta_{L,n}/\Gamma)_T. 
$$ 
This induces a corresponding map in cohomology, 
$$ 
(\pi_{e_{L/K}} \circ J_{L,K,n,T})^*: 
H^1(\cW(\bar\Delta_{L,n}/\Gamma)_T) \to 
H^1(\cW(\bar\Delta_{K,n}/\Gamma)_T). 
$$ 
To see this at the level of the filtrations of Proposition 
\ref{filtration}, notice that the map $J_{L,K,n}$ also induces a 
restriction map 
$$ 
r_{L,K,n}: {\rm C}(\cW(\bar\Delta_{L,n}/\Gamma),\Z) \to {\rm 
C}(\cW(\bar\Delta_{K,n}/\Gamma), \Z). 
$$ 
If we denote by $\cP^{(L)}_N$ and $\cP^{(K)}_j$ the respective 
filtrations, then we have restriction maps  
$r_{L,K,n}:\cP^{(L)}_{j e_{L/K}} \to \cP^{(K)}_j$. If we denote by 
$\delta_\ell$ the map $\delta_\ell(f)= f- f\circ T^\ell$, then the 
restriction also satisfies $r_{L,K,n} \circ \delta_{e_{L/K}} = 
\delta \circ r_{L,K,n}$, hence there is an induced map $\bar 
r_{L,K,n}: F^{(L,e_{L/K})}_j \to F^{(K)}_j$, where we have set 
$F^{(L,e_{L/K})}_j:= \cP^{(L)}_{j e_{L/K}}/\delta_{e_{L/K}} 
\cP^{(L)}_{(j-1) e_{L/K}}$. An argument analogous to the one used 
in the proof of Proposition \ref{filtration} shows that the 
$F^{(L,e_{L/K})}_j$ give a filtration of 
$$ H^1(\cW(\bar\Delta_{L,n}/\Gamma)_{T^{e_{L/K}}},\Z) 
=\varinjlim_j F^{(L,e_{L/K})}_j. $$ There is a corresponding map 
$${\rm C}(\cW(\bar\Delta_{L,n}/\Gamma),\Z)_T \to {\rm 
C}(\cW(\bar\Delta_{L,n}/\Gamma),\Z)_{T^{e_{L/K}}},$$ induced by 
the covering $\pi_{e_{L/K}}: 
\cW(\bar\Delta_{L,n}/\Gamma)_{T^{e_{L/K}}} \to 
\cW(\bar\Delta_{L,n}/\Gamma)_T$. On the level of filtrations this 
has the following description. The module $\cP^{(L)}_{j e_{L/K}}$ 
can be identified with the span of functions in $\cP^{(L)}_{j,s}$, 
$s=0,\ldots, e_{L/K}-1$, where $\cP^{(L)}_{j,s}= 
T^j(\cP^{(L)}_j)$. The inclusion $\iota: \cP^{(L)}_j 
\hookrightarrow \cP^{(L)}_{j e_{L/K}}$ as $\cP^{(L)}_{j,0}$ then 
satisfies $\delta\circ \iota=\iota\circ \delta_{e_{L/K}}$. We 
obtain an induced map $F^{(L)}_j \to F^{(L,e_{L/K})}_j$. Thus, the 
map induced by the covering $\pi_{e_{L/K}}$ also preserves the 
filtrations, and we obtain maps $F^{(L)}_j \to F^{(K)}_j$ that 
induce $(\pi_{e_{L/K}} \circ J_{L,K,T})^*$ on the direct limits. 
 
\noindent $\diamond$

\subsection{Dynamical cohomology} 
 
Let $\Delta$ be a directed tree on which the Schottky group 
$\Gamma$ acts, with the same assumptions as in the previous 
paragraph. 
 
On the set of ends $\partial \Delta$ we consider a measure $d\mu$ 
defined by first introducing the distance function $d(v) := {\rm 
dist}(v,x_0)$, for $v\in \Delta^0$ and $x_0$ the base point in 
$\Delta^0$ with respect to which the structure of directed graph 
on $\Delta$ is determined. Then the measure on $\partial\Delta$ is 
defined by assigning its value on the clopen set $V(v)$, given by 
the ends of all paths in $\Delta$ starting at a vertex $v$, to be 
$$ \mu(V(v))=q^{-d(v)-1}, $$ 
with $q={\rm card}(\O/m)$. 
 
\begin{prop} 
This induces a measure on $\cW(\Delta/\Gamma)$, with respect to 
which the shift map $T$ is measure preserving. 
\end{prop} 
 
\noindent\underline{Proof.} Notice that, if we identify the points 
of $V(v)$ with infinite paths starting at $v$, 
$$ V(v)=\{ w_0 w_1 \ldots w_n \ldots : s(w_0)=v, w_k\in (\Delta)^1_+ 
\}, $$ then the image $T(V(v))$ will have measure 
$\mu(T(V(v)))=\mu(V(v))/q$. In the case of walks starting at a 
vertex $v$, the map $T$ scales the measure of the set of walks 
starting with an edge $w\in (\Delta)^1_+$ by a factor $q^{-1}$ and 
the measure of the set of walks starting with an edge $w$ with 
$\bar w\in (\Delta)^1_+$ by a factor $q$. 
 
We can define a map from $\cW(\Delta/\Gamma)$ to $V\times V$, by 
splitting each doubly infinite sequence 
$$ \ldots w_{-n} w_{-n+1} \ldots w_{-1} w_0 w_1 \ldots w_\ell w_{\ell 
+1} \ldots $$ into the two sequences 
\begin{equation}\label{cutseq} 
 (w_0 w_1 \ldots w_\ell w_{\ell 
+1} \ldots, \,  \bar w_{-1} \bar w_{-2} \ldots \bar w_{-n+1} \bar 
w_{-n}\ldots ), 
\end{equation} 
each of which defines a point in the fundamental domain $V$, if we 
identify admissible sequences of edges in the fundamental domain 
$\Delta/\Gamma$ with admissible sequences of edges in $\Delta$ 
with the condition that the first edge $w_0$ (or $\bar w_{-1}$) 
lies in a chosen fundamental domain $\bar V$ of the action of 
$\Gamma$ which contains the base point $x_0$. Then the action of T 
maps \eqref{cutseq} to 
\begin{equation}\label{cutseqT} 
 (w_1 \ldots w_\ell w_{\ell 
+1} \ldots, \, \bar w_0  \bar w_{-1} \bar w_{-2} \ldots \bar 
w_{-n+1} \bar w_{-n}\ldots ), 
\end{equation} 
hence it scales the measure on one factor by $q$ and on the other 
factor by $q^{-1}$, so that the measure induced on 
$\cW(\Delta/\Gamma)$, by restricting to $V\times V$ the product 
measure on $\partial \Delta\times \partial \Delta$, is preserved 
by $T$. 
 
\noindent $\diamond$ 
 
Consider the free $\Z$-modules $\cP_n$, introduced in the proof of 
Proposition \ref{filtration}, of functions of at most $n+1$ future 
coordinates. We can realize $\cP_n\otimes_\Z \C$ as a vector 
subspace of $L^2(\partial\Delta,d\mu)$. The operator $\delta$ is 
bounded in norm, hence the $F_n = \cP_n /\delta \cP_{n-1}$ have 
induced norms and bounded inclusions $F_n\otimes \C 
\hookrightarrow F_{n+1}\otimes \C$, where the $F_n$ are the 
torsion free $\Z$-modules of Proposition \ref{filtration}. 
 
Using the inner product induced from $L^2(\partial \Delta,d\mu)$, 
we can split $F_n \cong F_{n-1} \oplus (F_{n-1}^\perp \cap F_n)$. 
 
\begin{defn}\label{def-dyncoh} 
The dynamical cohomology $\cH^1_{dyn}(\Delta/\Gamma)$ is the norm 
completion of the graded complex vector space 
$$ H^1_{dyn}(\Delta/\Gamma) := \oplus_{k\geq 0} Gr_k $$ 
with $Gr_n :=(F_n \cap F_{n-1}^\perp)$. 
\end{defn} 
 
We now construct a representation, by linear bounded operators, of 
a graph ${\rm C}^*$-algebra on this Hilbert space. 
 
\begin{prop}\label{repres} 
There is a representation of the algebra ${\rm 
C}^*(\Delta/\Gamma)$ by bounded linear operators on the Hilbert 
space $\cH^1_{dyn}(\Delta/\Gamma)$. 
\end{prop} 
 
\noindent\underline{Proof.} For $w\in (\Delta/\Gamma)^1_+$, we 
define linear operators $T_w$ of the form 
$$ (T_w f) (w_0 w_1 w_2 \ldots) = \left\{ \begin{array}{lr} f(w w_0 w_1 
w_2 \ldots) & r(w)=s(w_0) \\ 
0 & r(w)\neq s(w_0). \end{array} \right. $$ 
For $v\in 
\Delta/\Gamma^0$, let $P_v$ denote the projection 
$$ (P_v f)(w_0 w_1 w_2 \ldots) = \left\{ \begin{array}{lr} f(w_0 w_1 w_2 
\ldots) & s(w_0)=v \\ 
0 & s(w_0)\neq v, \end{array}\right. $$ 
We also define projections 
$\Pi_w$ of the form 
$$ (\Pi_w f)(w_0 w_1 w_2 \ldots) = \left\{ \begin{array}{lr} f(w_0 w_1 w_2 
\ldots) & w_0 = w \\ 
0 & w_0 \neq w. \end{array} \right. $$ Finally, we define linear 
operators $s_w$ and $S_w$ 
\begin{equation}\label{Sw} s_w := \sum_{w'} A(w,w') T_w 
\Pi_{w'} \ \ \ \text{ and } \ \ \ S_w:= \sqrt q\, s_w. 
\end{equation} 
These define bounded linear operators on $L^2(\partial\Delta, 
d\mu)$. The $S_w$ are partial isometries, satisfying $S_w S_w^* 
=\Pi_{w}$. Thus, the operators $\{ P_v, S_w \}$ form a 
Cuntz--Krieger family, since we obtain $P_v = \sum_{s(w)=v} \Pi_w 
= \sum_{s(w)=v} S_w S_w^*$ and $S_w^* S_w = \sum A(w,w') S_{w'} 
S_{w'}^*= \sum A(w,w') P_{w'} = P_{r(w)}$. 
 
The same argument given in \cite{CM}, Proposition 4.19, shows that 
$S_w \delta = \delta S_w$, so that they descend to bounded 
operators on $\cH^1_{dyn}(\Delta/\Gamma)$. 
 
\noindent $\diamond$ 
 
\subsection{Embedding of cohomologies} 
 
Let $\Delta$ be any of the trees $\bar\Delta_{K,n}$, $n\geq 0$, or 
$\Delta_\Gamma$. We construct a family of embeddings of the 
cohomology of the dual graph $\Delta_\Gamma/\Gamma$ in the 
dynamical cohomology $\cH^1_{dyn}(\Delta/\Gamma)$. Let 
$|\Delta_\Gamma|$ denote the geometric realization of the dual 
graph $\Delta_\Gamma/\Gamma$. 
 
\begin{thm} 
For each $N\geq 0$ there are embeddings $\phi_N$ of the cohomology 
of the dual graph $H^1(|\Delta_\Gamma/\Gamma|,\C)$ into the 
dynamical cohomology $\cH^1_{dyn}(\Delta/\Gamma)$. 
\end{thm} 
 
\noindent\underline{Proof.} We define maps $\phi_N$ from the 
cohomology $H^1(|\Delta_\Gamma/\Gamma|,\C)$ to 
$\cH^1_{dyn}(\Delta_K/\Gamma)$ in the following way: let $\{ 
\gamma_i \}_{i=1}^g$ be a chosen set of generators of the Schottky 
group $\Gamma$, $[\gamma_i]$ the corresponding homology classes in 
$H_1(|\Delta_\Gamma/\Gamma|,\Z)=\Gamma /[\Gamma,\Gamma]$, and 
$\eta_i$ the dual generators in cohomology. Let $w_i$ be the 
finite admissible word in the edges of the directed graph $\bar V$ 
(the fundamental domain for the $\Gamma$-action on $\Delta$) that 
represents the generator $\gamma_i$, with length $\ell_i = | w_i 
|$. Here, in the case of $\Delta=\bar\Delta_{K,n}$, we first 
notice that the first cohomology group of the dual graph is the 
same as $H^1(|\Delta_\Gamma '/\Gamma|,\Z)$, since the insertion of 
extra vertices does not change the topology of the graph. Since we 
have $\Delta_\Gamma'\subset \Delta_{K,n}$ we obtain the $w_i$ as 
above, as edges in the corresponding fundamental domain $\bar V$ 
for $\Delta_{K,n}$. 
 
We then set 
\begin{equation}\label{PhiN} 
 \phi_N (\eta_i)= P^\perp_{N \ell_i}\, \, 
\chi_{i,N}, 
\end{equation} 
where 
\begin{equation}\label{chi-in} 
\chi_{i,N}:=\chi_{\cW^+(\Delta/\Gamma,\underbrace{w_i\cdots 
w_i}_{N-times})} 
\end{equation} 
is the characteristic function of the set 
$\cW^+(\Delta/\Gamma,w_i\cdots w_i)$ of walks in $\Delta/\Gamma$ 
that begin with the word $w_i$  repeated $N$-times. The elements 
$\chi_{i,N}$ lie in $F_{N \ell_i}$. We denote by $P_k^\perp$ the 
orthogonal projection of $F_k$ onto $Gr_k$. 
 
The elements $\phi_N (\eta_i)$ are all linearly independent in 
$H^1_{dyn}(\Delta/\Gamma)$, hence the $\phi_N$ give linear 
embeddings of $H^1(|\Delta_\Gamma/\Gamma|,\C)$ into 
$\cH^1_{dyn}(\Delta/\Gamma)$. 
 
\noindent $\diamond$ 
 
\begin{cor} 
For any finite set of distinct $\{ N_k \}$, the set $\cup_k \{ 
\phi_{N_k}(\eta_i) \}_{i=1}^g$ consists of linearly independent 
vectors in $H^1_{dyn}(\Delta/\Gamma)$. Thus, we obtain an 
embedding 
\begin{equation}\label{Phiembedding} 
 \Phi=\oplus_N \phi_N : \bigoplus_N H^1(|\Delta_\Gamma/\Gamma|,\C) 
\to \cH^1_{dyn}(\Delta/\Gamma). 
\end{equation} 
\end{cor} 
 
Notice that the map $\Phi$ of \eqref{Phiembedding} does not 
preserve the graded pieces, due to the rescaling of the degrees by 
$\ell_i$ in \eqref{PhiN}. This has to be taken into account if we 
want to recover arithmetic information such as the local 
$L$--factors of \cite{Den} from the dynamical cohomology. For this 
reason, it may be necessary to blow up of some double points on 
the special fiber. 
 
\begin{lem}\label{blowup} 
It is always possible to reduce to the case where all the 
$\ell_i=\ell$, after blowing up a certain number of double points 
on the special fiber. 
\end{lem} 
 
\noindent\underline{Proof.}  Blowing-up a double point on the 
special fiber $C(k)$, for $k=\O/\m$, corresponds to introducing 
one extra vertex in the dual graph $\Delta_\Gamma/\Gamma$. This 
changes by one the lengths $\ell_i$ for those generators of 
$\Gamma$ for which the corresponding chain of edges $w_i$ passes 
through the newly inserted vertex. 
 
\noindent $\diamond$ 
 
Thus, possibly after blowing up some double points, we obtain: 
$$\phi_N : H^1(|\Delta_\Gamma/\Gamma|,\C) \hookrightarrow 
Gr_{N\ell} \subset \cH^1_{dyn}(\Delta_K/\Gamma), $$ for some 
$\ell\geq 1$. 
 
\subsection{Spectral triples} 
 
In this paragraph we show that we can associate to a given Mumford 
curve a family of spectral triples $(\cA_n,\cH_n,D_n)$, for $n\geq 
-1$. Each triple in this family corresponds to the choice of a 
graph $\Delta =\bar\Delta_{K,n}$, for $n\geq 0$ and 
$\Delta=\Delta_\Gamma$ for $n=-1$. The Dirac operator in these 
spectral triples depends only on the graded structure of the space 
$H^1_{dyn}(\Delta/\Gamma)$. 
 
\medskip 
 
Recall that a {\em spectral triple} $(\cA,\cH,D)$ is a triple of a 
${\rm C}^*$-algebra, with a representation in the algebra of 
bounded operators on the Hilbert space $\cH$, and a Dirac operator 
$D$, which is a self-adjoint unbounded operator acting on $\cH$, 
such that $(\lambda-D)^{-1}$ is a compact operator for all 
$\lambda\notin \R$ and the commutators $[a,D]$ are bounded 
operators for all $a\in \cA_0$, with $\cA_0$ a dense involutive 
subalgebra of $\cA$. 
 
Recall also that a spectral triple  $(\cA,\cH,D)$ has an 
associated family of zeta functions of the form 
\begin{equation}\label{zeta} 
\zeta_{a,|D|}(z)=Tr(a|D|^{-z})=\sum_{\lambda\in \Sp(|D|)\setminus 
\{ 0 \} } Tr(a \Pi_\lambda) \lambda^{-z}, 
\end{equation} 
with $a\in \cA_0 \cup [D,\cA_0]$, and with $\Pi_\lambda$ the 
spectral projection on $\lambda\in \Sp(|D|)$. The properties of 
these zeta functions are related to the notion of dimension 
spectrum for the spectral triple (the set of poles of the 
$\zeta_{a,|D|}$) and to the local index formula of Connes and 
Moscovici. For the purpose of this paper, we will extend the zeta 
functions \eqref{zeta} to the case where the element $a$ is a weak 
limit of certain sequences of elements in $\cA$, since, in our 
construction, this is the type of zeta functions that recovers the 
arithmetic invariant given by the local Euler factor. 
 
There are corresponding two-variable zeta functions 
\begin{equation}\label{zeta-s} 
\zeta_{a,|D|}(s,z)=\sum_{\lambda\in \Sp(|D|)} Tr(a \Pi_\lambda) 
(s+\lambda)^{-z}, 
\end{equation} 
and associated regularized determinants 
\begin{equation}\label{reg-det} 
\det_{\infty, a,|D|} (s) = \exp \left( - \frac{d}{dz} 
\zeta_{a,|D|} (s,z)|_{z=0} \right). 
\end{equation} 

\medskip 
 
We consider the Hilbert space 
\begin{equation}\label{H-sp} 
\cH=\cH^1_{dyn}(\Delta/\Gamma)\oplus \cH^1_{dyn}(\Delta/\Gamma). 
\end{equation} 
On this space we consider the diagonal action of ${\rm 
C}^*(\Delta/\Gamma)$. We also introduce the notation $Gr_{n,-}:= 
Gr_n\oplus 0$ and $Gr_{n,+}:= 0 \oplus Gr_n$. 
 
We define the Dirac operator $D$ acting on $\cH$ by setting 
\begin{equation}\label{D} 
 D|_{Gr_{+,n}}=n \ \ \ \ \  D|_{Gr_{-,n}}=-n-1. 
\end{equation} 
 
\begin{prop}\label{triple1} 
The data $({\rm C}^*(\Delta/\Gamma), \cH^1_{dyn}(\Delta/\Gamma) 
\oplus \cH^1_{dyn}(\Delta/\Gamma), D)$ determine a spectral triple 
in the sense of Connes \cite{Connes}. 
\end{prop} 
 
\noindent\underline{Proof.}  In order to obtain a spectral triple 
we need to check the compatibility requirement between the Dirac 
operator $D$ and the action of the algebra ${\rm 
C}^*(\Delta/\Gamma)$. It is sufficient to check that the 
commutators $[D,S_w]$ are bounded operators. This follows easily 
since $S_w : Gr_{\pm, n} \to Gr_{\pm, n-1}$, so that $[ D, S_w ]f= 
\mp S_w f$. The remaining properties are easily verified. 
 
\noindent $\diamond$ 
 
We can modify slightly the above construction, in order to take 
into account the scaling factor $\ell$ in the grading between 
$\oplus_N H^1(|\Delta_\Gamma/\Gamma|,\C)$ and its image under 
$\Phi$ in $\cH^1_{dyn}(\Delta/\Gamma)$. 
 
\begin{cor}\label{D-nl} 
We modify the operator $D$ of \eqref{D}, by setting 
\begin{equation}\label{D-2} 
D|_{Gr_{+,n}} =\frac{n}{\ell}\frac{2\pi}{\log
q}  \ \ \ \ \ D|_{Gr_{-,n}} =\frac{-(n+1)}{\ell}\frac{2\pi}{\log q}.
\end{equation} 
Here $q$ is the cardinality of the residue field 
$k(v)$ and $\ell$ is the length of all the words representing the 
generators of $\Gamma$, possibly after blowing up some points on 
the special fiber. With this modified operator $D$, we still 
obtain a spectral triple $$({\rm C}^*(\Delta/\Gamma), 
\cH^1_{dyn}(\Delta/\Gamma) \oplus \cH^1_{dyn}(\Delta/\Gamma), D).$$ 
\end{cor} 
 
\subsection{Local $L$-factor} 
 
Let $X$ be a curve over a global field $\mK$. We assume 
semi-stability at all places of bad reduction. The local Euler 
factor at a place $v$ has the following description 
(\cite{Serre}): 
\begin{equation}\label{L-factor} 
 L_v (H^1(X),s)= \det\left( 1-Fr_v^* N(v)^{-s} | H^1(\bar X, 
\Q_\ell)^{I_v} \right)^{-1}. 
\end{equation} 
Here $Fr_v^*$ is the geometric Frobenius acting on $\ell$-adic 
cohomology of $\bar X=X \otimes \Sp(\bar\mK)$, with $\bar\mK$ the 
algebraic closure and $\ell$ a prime with $(\ell,q)=1$, where $q$ 
is the cardinality of the residue field $k(v)$ at $v$. We 
denote by $N$ the norm map. The determinant is evaluated on the 
inertia invariants $H^1(\bar X, \Q_\ell)^{I_v}$ at $v$ (all of 
$H^1(\bar X, \Q_\ell)$ when $v$ is a place of good reduction). 
 
Suppose $v$ is a place of $k(v)$-split degenerate reduction. Then 
the completion of $X$ at $v$ is a Mumford curve $X_\Gamma$. In 
this case, the Euler factor \eqref{L-factor} takes the following 
form: 
\begin{equation}\label{loc-factor} 
L_v(H^1(X_\Gamma),s)= \prod_\lambda (1-\lambda q^{-s})^{-\dim 
H^1(X_\Gamma)_\lambda^{I_v}}=(1-q^{-s})^{-g}, 
\end{equation} 
since the eigenvalues $\{ \lambda \}$ of the Frobenius, in this 
case, are all $\lambda=1$. Here we denote by 
$H^1(X_\Gamma)_\lambda^{I_v}$ the eigenspaces of the Frobenius. 
 
Deninger in \cite{Den} and \cite{Den2} obtained the local factor
\eqref{L-factor} as a regularized determinant over an infinite
dimensional cohomological theory. 
 
In the case of Mumford curves, Deninger's calculation can be 
recast in terms of the data of the spectral triples of Proposition 
\ref{triple1} and of the embeddings $\Phi^\pm$ of cohomologies.

We consider the operator $iD$ (an imaginary rotation of the Dirac
operator) and the zeta functions
\begin{equation}\label{zetapm}
\begin{array}{ll}
\zeta_{a,iD,+} (s,z) := & \sum_{\lambda \in \Sp(iD)\cap i[0,\infty)}
\Tr(a\Pi_\lambda) (s+\lambda)^{-z} \\[3mm]
\zeta_{a,iD,-} (s,z) := & \sum_{\lambda \in \Sp(iD)\cap i(-\infty,0)}
\Tr(a\Pi_\lambda) (s+\lambda)^{-z}.
\end{array}\end{equation}
Then we have a regularized determinant
\begin{equation}\label{det-pm}
\det_{\infty,a,iD}(s) := \exp\left( -\zeta^\prime_{a,iD,+}(s,0) \right) 
\exp\left(-\zeta^\prime_{a,iD,-}(s,0) \right)
\end{equation}

\begin{thm}\label{L-factor-1} 
Let $\pi(\cV)$ be the orthogonal projection of 
$\cH^1_{dyn}(\Delta/\Gamma)\oplus \cH^1_{dyn}(\Delta/\Gamma)$ onto 
the graded subspace $\cV= {\rm Im}(\Phi^-) \oplus {\rm 
Im}(\Phi^+)$, where we denote by $\Phi^\pm$ the maps $\Phi\oplus 
0$ and $0\oplus\Phi$. Then the regularized determinant 
\eqref{det-pm}, with $a=\pi(\cV)$ and $D$ the Dirac operator of 
Corollary \ref{D-nl} satisfies 
\begin{equation}\label{det-pm-L} 
\det_{\infty,\pi(\cV),iD}( s ) = 
L_v (H^1(X_\Gamma),s)^{-1}. 
\end{equation} 
\end{thm} 
 
\noindent\underline{Proof.} 
When we compute the zeta functions \eqref{zetapm} for the Dirac operator 
of Corollary \ref{D-nl}, and $a=\pi(\cV)$, we obtain
$$ \begin{array}{ll}
 \zeta_{\pi(\cV),iD,+} (s,z) = & \sum_{n=0}^\infty \Tr
(\pi(\cV)\Pi_{+,n\ell} )  (\gamma (\tau + n))^{-z} \\[3mm]
 \zeta_{\pi(\cV),iD,-} (s,z) = & \sum_{n=0}^\infty \Tr
(\pi(\cV)\Pi_{-,n\ell}) (\gamma (\tau - n))^{-z} - \Tr
(\pi(\cV)\Pi_{-,0}) (\tau \gamma)^{-z},  
\end{array} $$
for $\gamma = \frac{2\pi i}{\log q}$ and $\tau = \frac{ \log q}{2\pi
i} s$, with choice of arguments $-\pi < \arg \gamma (\tau + n) < \pi$,
as in \cite{Den2}. 

Furthermore, we have $Tr(\pi(\cV)\Pi_{\pm,n\ell})= \dim
(Gr_{\pm,n\ell} \cap \cV)= g$.  
In fact, the space $Gr_{+,n\ell} \cap \cV={\rm Im}(\Phi^+)$ is 
generated by $0\oplus \chi_{i,n}$ for $i=1,\ldots,g$ and 
$Gr_{-,n\ell} \cap \cV={\rm Im}(\Phi^-)$ is generated by the 
element $\chi_{i,n}\oplus 0$ for $i=1,\ldots,g$. 

The result then follows the calculation of the regularized 
determinant given in \cite{Den2}. Namely, we obtain
$$ (1-q^{-s})^{-g} = (\tau\gamma)^{-g} \exp(-g \zeta_\gamma'(\tau,0))
\exp(-g \zeta_{-\gamma}'(-\tau,0)), $$
which is exactly the regularized determinant $\det_\infty
(s-\Theta_q)$ computed in \cite{Den2} for the spectrum (with
multiplicity $g$)
$$
\Sp(s-\Theta_q)=\left\{ \frac{2\pi i}{\log q} \left( \frac{s\log
q}{2\pi i} + n \right) : n\in \Z \right\}. 
$$

\noindent $\diamond$

It is very interesting to notice an important difference between the
archimedean and non-archimedean cases. At the archimedean prime (\cf
\cite{CM}, \cite{Den1}) the local factor is described in terms of zeta
functions for a Dirac operator $D$. On the other hand, at the
non-archimedean places, in order to get the correct normalization as
in \cite{Den2}, we need to introduce a rotation of the Dirac operator
by the imaginary unit, $D\mapsto iD$. This rotation corresponds to 
the Wick rotation that moves poles on the real line to
poles on the imaginary line (zeroes for the local factor) and appears 
to be a manifestation of a rotation from Minkowskian to Euclidean
signature $it \mapsto t$, as already remarked by Manin 
(\cite{Man-zeta} p.135), who wrote that {\em ``imaginary time motion'' 
may be held responsible for the fact that zeroes of $\Gamma(s)^{-1}$ are
purely real whereas the zeroes of all non-archimedean Euler factors
are purely imaginary}. This seems to hint to the existence of a 
more refined construction involving Minkowskian geometry, where  
the rotation $D\mapsto iD$ could be interpreted as a rotation 
$it \mapsto t$ of an infinitesimal length element $D^{-1} \sim ic\, dt$. 
A more precise treatment would require adapting the structure of spectral 
triple to the case of Minkowskian signature. Another piece
of supporting evidence for the idea that a more refined construction
should involve Minkowskian geometry comes from the cohomological
construction of \cite{KC}. In fact, in \cite{CM}, we only used part
of the full symmetry group determined by the Lefschetz module
structure, namely the part corresponding to real hyperbolic
geometry, so as to match the results of \cite{Man-3d}. The full symmetry
group is $\SL(2,\R)\times \SL(2,\R)$, which is the isometry group not
of 3-dimensional real hyperbolic geometry, but of its Minkowskian
version, anti de Sitter space ${\rm AdS}_{2+1}$. On the relation
between the results of \cite{Man-3d} and ${\rm AdS}_{2+1}$ see also
\cite{ManMar}.

\subsection{AF core} 

In order to understand the zeta functions 
$\zeta_{\pi(\cV),D,\pm}(\tau,z)$ of Theorem \ref{L-factor-1} in 
terms of the spectral triple, we still need to express the 
operator $\pi(\cV)$ in terms of elements in the algebra ${\rm 
C}^*(\Delta/\Gamma)$. 
 
The graph ${\rm C}^*$-algebra ${\rm 
C}^*(\Delta/\Gamma)$ contains an AF core given by the AF algebra
obtained as the norm closure of $\cup_n {\mathcal 
F}_n$, where the finite dimensional algebras ${\mathcal F}_n$ are
given by  
\begin{equation}\label{AF-core} 
{\mathcal F}_n= {\rm span}\{ S_\mu S_\nu^* : \mu,\nu \in 
\cP^n(\Delta/\Gamma), t(\mu)=t(\nu) \}. 
\end{equation} 
Here we used the notation $S_\mu=S_{w_1}\cdots S_{w_k}$, for 
$\mu=w_1\cdots w_k$. The AF core can be identified with the fixed 
point algebra of the gauge action (\cf \cite{BPRS}), 
\begin{equation}\label{AF-core-cl} 
\overline{\bigcup_n {\mathcal F}_n} \cong {\rm 
C}^*(\Delta/\Gamma)^{U(1)}. 
\end{equation} 
 
Let $\{ w_i \}_{i=1}^g$ be the words corresponding to the 
generators of $\Gamma$, which we assume all of equal length 
$\ell$, possibly after some blow ups. Let $S_{w_i}$ be the 
corresponding operators in ${\rm C}^*(\Delta/\Gamma)$. The 
operators $S_{w_i}^n {S_{w_i}^*}^n$ belong to the subalgebra 
${\mathcal F}_{n\ell}$ in the AF core of ${\rm 
C}^*(\Delta/\Gamma)$. 

Each $Q_{i,n}= S_{w_i}^n {S_{w_i}^*}^n$ acts 
on $L^2(\partial\Delta,d\mu)$ as multiplication by the 
characteristic function $\chi_{i,n}$ of \eqref{chi-in}, hence
$Q_{i,n}$ maps $Gr_{\pm,n\ell}$ to itself, with range the one dimensional
subspace of $Gr_{\pm,n\ell}$ spanned by $\chi_{i,n}$. Thus, the 
operator $Q_n=\sum_i Q_{i,n}$ projects $Gr_{\pm,n\ell}$ onto 
the $g$-dimensional subspace $Gr_{\pm,n\ell} \cap \cV$.

For the Dirac operator of \eqref{D-2}, we write
$D=\sum_{n\geq 0} \Pi_n \lambda_n$, where
with $\lambda_n =2\pi n/\ell \log q$ and 
$\Pi_n = \Pi_{+,n}\oplus \Pi_{-,n-1}$. We then have the 
following result. 

\begin{prop}\label{L-factor-2}
the zeta function 
$$ \zeta_{\pi(\cV),|D|}(z) = Tr(\pi(\cV) |D|^{-z}) $$ 
can be written in the form
\begin{equation}\label{L-2}
\zeta_{\pi(\cV),|D|}(z) =\Tr\left( \sum_{n>0} Q_{i,n} \Pi_{n\ell} 
\lambda_n^{-z} \right) 
\end{equation} 
with the $Q_{i,n}$ in the AF core of ${\rm C}^*(\Delta/\Gamma)$.
\end{prop}

\section{Foam spaces}

We now consider the local factor \eqref{L-factor} in the more general
case, where we drop the assumption that  $v$ is a place of
$k(v)$-split degenerate reduction. In this case, we no longer have a
p--adic uniformization of the completion of $X$ at $v$ as a Mumford
curve, and correspondently, the local factor is no longer determined
solely in terms of the combinatorics of the dual graph, but it depends
essentially on extra geometric information on the nature of the
degeneration. In particular, the inertia invariants $H^1(\bar X, 
\Q_\ell)^{I_v}$ are described only partly by the cohomology of the
dual graph, with the extra information provided by the cohomology
of the single components of the dual fiber, which in the general case
will no longer be just rational curves.

More precisely, if we denote by $H^1(\bar X)_\lambda^{I_v}$ the
eigenspace of the geometric Frobenius, with eigenvalue $\lambda$, we
can write the Euler factor in the form 
\begin{equation}\label{loc-factor-lambda} 
L_v(H^1(X),s)= \prod_\lambda (1-\lambda q^{-s})^{-d_\lambda}, 
\end{equation} 
with $d_\lambda = \dim H^1(\bar X)_\lambda^{I_v}$.
Deninger's description of the local factor as regularized determinant
holds in this more general case, in the form $\det_\infty (s -\Theta)$
where $\Theta$ has spectrum $\{ \alpha_\lambda +\frac{2\pi i n}{\log
q} \}$, with $n\in \Z$, $\lambda \in
\Sp(Fr_v^*)$, and $q^{\alpha_\lambda}=\lambda$.

We want to modify the graphs $\Delta/\Gamma$ considered in the
previous sections, in such a way that the corresponding dynamical
cohomology will contain a linear subspace isomorphic to an infinite
direct sum $\oplus_N H^1(\bar X, 
\Q_\ell)^{I_v}$, and such that the construction of the spectral triple
and the derivation of the regularized determinant described in the
case of Mumford curves will extend to this case to recover
\eqref{loc-factor-lambda}. 

Using the exact sequence of \cite{Mor} p.110-111, we obtain an
identification 
\begin{equation}\label{cohom1}
H^1(\bar X, \Q_\ell)^{I_v}\cong H^1(|\Delta_\Gamma
/\Gamma|)\otimes_{\sigma_\ell} \Q_\ell \oplus H^1(X^{[0]}_v)
\otimes_{\sigma_\ell} \Q_\ell, 
\end{equation}
where $\sigma_\ell : \Q_\ell \to \C$ is a fixed embedding of $\Q_\ell$
in $\C$, for a prime $\ell$ with $(\ell, q)=1$, and we denote by 
$X^{[0]}_v$ the disjoint union of the components of the special fiber. 
In the case of $k(p)$-split reduction, where all components are 
$\P^1$'s, \eqref{cohom1} is simply identified with $H^1(|\Delta_\Gamma
/\Gamma|)\otimes_{\sigma_\ell} \Q_\ell$ as in the previous sections.
The finite decomposition $H^1(\bar X, \Q_\ell)^{I_v} =\oplus_\lambda
H^1(\bar X, \Q_\ell)^{I_v}_\lambda$ in eigenvalues of the geometric
Frobenius provides corresponding spaces $H^1(|\Delta_\Gamma
/\Gamma|)_\lambda$ and $H^1(X^{[0]}_v)_\lambda$ of dimensions
$d^\Gamma_\lambda$ and $d^0_\lambda$, respectively, with
$d^\Gamma_\lambda + d^0_\lambda = d_\lambda$. 

We choose vertices $x_\lambda$ (not necessarily distinct) of
$\Delta_\Gamma /\Gamma$, and attach to the vertex $x_\lambda$
new outgoing edges $w_{i,\lambda}^0$, with $i=1,\ldots
d_\lambda^0$. We denote by $E_v$ the oriented 
graph obtained via this construction, after appending tails to all
sinks. 

\begin{rem}{\em
For certain classes of examples, our graph $E_v$ can be embedded as a
subgraph of the ``foam space'' defined in \cite{CMZ}. The foam space
is a graph $F_v$ associated to the fiber $X_v$ of an arithmetic surface
$\mX$ over $\Sp(\O_{\mK})$, obtained by replacing the special
fiber $X_v$ with an infinite series of blowups of its $\mF_q$-points.
The graph $F_v$ is the limit of the dual graphs associated to
this series of blow-ups (\cf \cite{ManC} \S 35). For this reason, we
think of the graphs $E_v$ as a generalization of ``foam spaces''.
}\end{rem}
 
To our foam space $E_v$, we associate the corresponding dynamical
cohomology $\cH^1(E_v)$ as in the previous sections, and the graph
${\rm C}^*$-algebra ${\rm C}^*(E_v)$. The argument of Proposition
\ref{triple1} extends to this case and gives a spectral triple
$$ ({\rm C}^*(E_v), \cH^1(E_v)\oplus \cH^1(E_v), D). $$

We now define embeddings of cohomology groups as follows. 
Let $\omega_{i,\lambda}$, for $i=1,\ldots,d^\Gamma_\lambda$ be loops
of edges in $\Delta_\Gamma/\Gamma$, with
$|\omega_{i,\lambda}|=\ell_{i,\lambda}$, representing homology classes
dual to a basis $\{ \eta_{i,\lambda}^\Gamma \}$ of
$H^1(|\Delta_\Gamma/\Gamma|)_\lambda$. Up to adding vertices to the
graph $\Delta_\Gamma/\Gamma$ by blowing up double points in the closed
fiber, we can assume that all the $\ell_{i,\lambda}=\ell$. Adding
vertices in this way does not change $H^1(X^{[0]}_v)$, since the 
components of the closed fiber that correspond to the new vertices
all have trivial $H^1$.

We consider then the linear embedding 
$$ \Phi_{N,\lambda}^\Gamma: H^1(|\Delta_\Gamma/\Gamma|)_\lambda
\hookrightarrow Gr_{N \ell} \subset \cH^1(E_v), $$
given by 
$$ \Phi_{N,\lambda}^\Gamma(\eta_{i,\lambda}^\Gamma) =
P_{N\ell_{\lambda}}^\perp 
\chi_{\cW^+(E_v, \underbrace{\omega_{i,\lambda} \cdots
\omega_{i,\lambda}}_{N-times})}. $$
We also consider the embeddings 
$$ \Phi_{N,\lambda}^0 : H^1(X^{[0]}_v)_\lambda \hookrightarrow Gr_{N\ell}
\subset \cH^1(E_v), $$ 
$$ \Phi_{N,\lambda}^0(\eta_{i,\lambda}^0) = P_{N}^\perp 
\chi_{\cW^+(E_v, \underbrace{w_{i,\lambda}^0 \cdots
w_{i,\lambda}^0}_{N\ell-times})}, $$
where the $\eta_{i,\lambda}^0$ form a basis of
$H^1(X^{[0]}_v)_\lambda$ and the $w_{i,\lambda}^0$ are the
corresponding oriented edges of $E_v$. 
Let $\Phi_\lambda^\Gamma = \oplus_N \Phi_{N,\lambda}^\Gamma$ and
$\Phi_\lambda^0 = \oplus_N \Phi_{N,\lambda}^0$, and let $\Phi_\lambda
= \Phi_\lambda^\Gamma \oplus \Phi_\lambda^0$. With $\Phi_\lambda^\pm$
defined as the $\Phi^\pm$ in Theorem \ref{L-factor-1}, we denote by
$\cV_\lambda = {\rm Im}(\Phi_\lambda^-) \oplus {\rm 
Im}(\Phi_\lambda^+)$, and by $\pi(\cV_\lambda)$ the corresponding
orthogonal projection.

We then extend the result of Theorem \ref{L-factor-1} to this more
general setting. 
   
\begin{thm} \label{L-factor-3}
Consider the regularized determinants  
\eqref{det-pm}, with $a_\lambda=\pi(\cV_\lambda)$ and $D$ the Dirac
operator of Corollary \ref{D-nl} for the spectral triple $({\rm
C}^*(E_v), \cH^1(E_v)\oplus \cH^1(E_v), D)$. We obtain  
\begin{equation}\label{det-lambda} 
\prod_\lambda \det_{\infty,\pi(\cV_\lambda),iD}(s) = L_v (H^1(X),s)^{-1}. 
\end{equation}
The operators $\pi(\cV_\lambda)$ are
related to the AF core of the ${\rm C}^*$-algebra
${\rm C}^*(E_v)$ as in \eqref{L-2}.
\end{thm} 

\noindent\underline{Proof.} 
We compute the zeta functions \eqref{zetapm} for the Dirac operator of
the spectral triple $({\rm C}^*(E_v), \cH^1(E_v)\oplus \cH^1(E_v),
D)$, with $a=\pi(\cV)$. We obtain 
$$ \begin{array}{ll}
\zeta_{\pi(\cV_\lambda),iD,+} (s,z) = & \sum_{n=0}^\infty \Tr
(\pi(\cV_\lambda)\Pi_{+,n\ell} )  
(\gamma (\tau_\lambda + n))^{-z} \\[3mm]
 \zeta_{\pi(\cV_\lambda),iD,-} (s,z) = & \sum_{n=0}^\infty \Tr
(\pi(\cV_\lambda)\Pi_{-,n\ell}) 
(\gamma (\tau_\lambda - n))^{-z} - \Tr (\pi(\cV_\lambda)\Pi_{-,0})
(\tau_\lambda \gamma)^{-z},  
\end{array} $$
for $\gamma = \frac{2\pi i}{\log q}$, $\tau_\lambda = \frac{\log
q}{2\pi i} (s -\alpha_\lambda)$, and $q^{\alpha_\lambda}=\lambda$, and
with choice of arguments as in \cite{Den2}. 

By construction, we have $Tr(\pi(\cV_\lambda)\Pi_{\pm,n\ell})= \dim
(Gr_{\pm,n\ell} \cap \cV_\lambda)= d_\lambda$, hence the left hand
side of \eqref{det-lambda} is 
the regularized determinant $det_\infty(s-\Theta_q)$ computed in
\cite{Den2}, with spectrum (with
multiplicities $d_\lambda$)
$$
\Sp(s-\Theta_q)=\left\{ \frac{2\pi i}{\log q} \left( \frac{\log
q}{2\pi i}(s- \alpha_\lambda) + n \right) : n\in \Z, \lambda \in
\Sp(Fr_v^*) \right\}.  
$$
The expression of the operators $\pi(\cV_\lambda)$ in terms 
of operators in the AF core of the ${\rm C}^*$-algebra
${\rm C}^*(E_v)$ is analogous to the case of Mumford
curves. 

\noindent $\diamond$

\bigskip
\bigskip
\bigskip

\noindent{\bf Acknowledgment.} Part of this work was done during 
visits of the first author to the Max Planck Institute in Bonn and 
of the second author to Florida State University and University of 
Toronto. We thank these institutions for hospitality and support. 
We are very grateful to Christopher Deninger for a crucial remark
about normalizations. 
This research has been partially supported by NSERC grant 72016789 
and by Humboldt Foundation and the German Government (Sofja 
Kovalevskaya Award).

\noindent {\bf Caterina Consani}, Department of Mathematics, 
University of Toronto, Canada. 
 
\noindent email: kc\@@math.toronto.edu 
 
\vskip .2in 
 
\noindent {\bf Matilde Marcolli}, Max--Planck--Institut f\"ur 
Mathematik, Bonn Germany. 
 
\noindent email: marcolli\@@mpim-bonn.mpg.de

\end{document}